\documentclass{article}

\usepackage[utf8]{inputenc}
\usepackage[T1]{fontenc}
\usepackage{fullpage}
\usepackage{xargs}
\usepackage{amsmath,amssymb,fancybox,pgf,array}

\newcommand{\virg}[1]{``#1''}

\renewcommand{\b}[1]{{\bf #1}}

\newcommand{\mmfunz}[5]{
$$ #1 : \left\{ \begin{array}{ccl}
 #2 & \rightarrow & #3 \\
 #4 & \mapsto& #5   \end{array} \right. $$}
\newcommand{\fz}[3]{#1:\, #2 \rightarrow #3}

\renewcommand{\r}[1]{(\ref{#1})}
\newcommand{\bi}{\begin{itemize}}
\newcommand{\ei}{\end{itemize}}
\newcommand{\be}{\begin{enumerate}}
\newcommand{\ee}{\end{enumerate}}

\newcommand{\bd}{\begin{description}}
\newcommand{\ed}{\end{description}}
\renewcommand{\i}{\item}

\newcommand{\bqn}{\begin{eqnarray}}
\newcommand{\eqn}{\end{eqnarray}}
\newcommand{\eqnn}{\nonumber\end{eqnarray}}
\newcommand{\eqnl}[1]{\label{#1}\end{eqnarray}}

\newcommand{\nn}{\nonumber\\}

\newcommand{\ba}[1]{\begin{array}{#1}}
\newcommand{\ea}{\end{array}}

\newcommand{\R}{\mathbb{R}}

\newcommand{\N}{\mathbb{N}}

\newcommand{\bproof}{\begin{proof}}
\newcommand{\eproof}{\end{proof}}
\newtheorem{Theorem}{\bf Theorem}
\newtheorem{lemma}[Theorem]{\bf Lemma}
\newtheorem{corollary}[Theorem]{\bf Corollary}
\newtheorem{definition}[Theorem]{\bf Definition}
\newtheorem{proposition}[Theorem]{\bf Proposition}
\newtheorem{remark}[Theorem]{\bf Remark}
\newenvironment{proof}{{\it Proof.}~~}{\hfill$\square$}

\newcommand{\bt}{\begin{Theorem}}
\newcommand{\et}{\end{Theorem}}
\newcommand{\bl}{\begin{lemma}}
\newcommand{\el}{\end{lemma}}
\newcommand{\bp}{\begin{proposition}}
\newcommand{\ep}{\end{proposition}}
\newcommand{\bc}{\begin{corollary}}
\newcommand{\ec}{\end{corollary}}
\newcommand{\bdeff}{\begin{definition}}
\newcommand{\edeff}{\end{definition}}
\newcommand{\brem}{\begin{remark}\rm}
\newcommand{\erem}{\end{remark}}
\newcommand{\auth}[1]{{\sc #1}}
\newcommand{\tit}[1]{{\rm #1}}
\newcommand{\titl}[1]{{\it #1}}
\newcommand{\jou}[1]{{\it #1}}

\newcommand{\pp}[1]{pp.~#1}

\newcommand{\lam}{\lambda}
\newcommand{\al}{\alpha}
\newcommand{\eps}{\varepsilon}
\newcommand{\ga}{\gamma}

\newcommand{\Id}{\mathrm{Id}}

\newcommand{\Pt}[1]{\left( #1 \right)}
\newcommand{\Pg}[1]{\left\{ #1 \right\}}
\newcommand{\Pq}[1]{\left[ #1 \right] }

\newcommand{\Pabs}[1]{{\Big \vert}  #1 {\Big \vert}}

\newcommand{\dt}{{\Delta t}}

\newcommand{\Mu}{\mathcal{M}}

\newcommand{\schema}[1]{\b{\sc #1}}

\newcommand{\Chi}{{\chi}}

\renewcommand{\dt}{{\Delta t}}
\newcommand{\weak}{\rightharpoonup}

\newcommand{\supp}{\mathrm{supp}}

\begin{document}
\newcommand{\sotto}{\leq}
\newcommand{\gwbase}{W^{a,b}_p}
\newcommandx{\gw}[4][1=a,2=b,3=p]{W^{#1,#2}_{#3}\Pt{#4}}

\newcommand{\Muu}{{\Mu_{0}^{ac}}}
\renewcommand{\P}{\mathcal{P}}

\title{Generalized Wasserstein distance and its application to transport equations with source}


\author{Benedetto Piccoli\thanks{Department of Mathematical Sciences, Rutgers University - Camden, Camden, NJ. {\tt piccoli@camden.rutgers.edu}}, Francesco Rossi\thanks{Aix-Marseille Univ, LSIS, 13013, Marseille, France. {\tt francesco.rossi@lsis.org}}}

\maketitle

\begin{abstract}
\noindent In this article, we generalize the Wasserstein distance to measures with different masses. We study the properties of such distance. In particular, we show that it metrizes weak convergence for tight sequences.\\

\noindent We use this generalized Wasserstein distance to study a transport equation with source, in which both the vector field and the source depend on the measure itself. We prove existence and uniqueness of the solution to the Cauchy problem when the vector field and the source are Lipschitzian with respect to the generalized Wasserstein distance.

\end{abstract}

\vspace{.3cm}


\vspace{1cm}

The problem of optimal transportation, also called  Monge-Kantorovich problem, has been intensively studied in the mathematical community. Related to such problem, the definition of the Wasserstein distance in the space of probability measure has revealed itself to be a powerful tool, in particular for dealing with dynamics of measures (like the transport PDE, see e.g. \cite{gradient}). For a complete introduction to Wasserstein distances, see \cite{old-new,villani}.

The main limit of this approach, at least for its application to dynamics of measures, is that the Wasserstein distance $W_p(\mu,\nu)$ is defined only if the two measures $\mu,\nu$ have the same mass. For this reason, in this article we first define a generalized Wasserstein distance $\gw{\mu,\nu}$, combining the standard Wasserstein and $L^1$ distances. In rough words, for $\gw{\mu,\nu}$ an infinitesimal mass $\delta\mu$ of $\mu$ (or $\nu$) can either be removed at cost $a |\delta\mu|$, or moved from $\mu$ to $\nu$ at cost $bW_p(\delta\mu,\delta\nu)$.


After proving simple but important properties of this distance, we show that it is the natural distance to state existence and uniqueness of the solution to the following dynamics of measures:
\bqn
\begin{cases}
\partial_t\mu+\nabla\cdot(v\Pq{\mu}\, \mu)=h\Pq{\mu},\\
\mu_{|_{t=0}}=\mu_0.
\end{cases}
\label{e-cauchy}
\eqn
Such equation is intensively used for modelling of crowd dynamics, where $\mu_t$ represents a pedestrian density (see e.g. \cite{ben-multiscale,evers,maury1,maury2,nostro,ben-ARMA,ben-domain,tosin}). Several authors have studied \r{e-cauchy} without source terms, i.e. $h\equiv 0$, showing that it is very convenient to use the standard Wasserstein distance in this framework. In particular, in \cite{ambrosio} the authors prove existence and uniqueness of the solution under Lipschitzianity of $v$ with respect to $\mu$, when the space of measures (of given fixed mass) is endowed with the Wasserstein distance. See also \cite{crippa}. We also showed that these hypotheses give convergence of numerical schemes \cite{nostro}.

The main limit of the approach based on the standard Wasserstein distance is that it cannot encompass the case of a source $h$. Indeed, in this case the mass of the measure $\mu_t$ varies in time, hence in general $W_p(\mu_t,\mu_s)$ is not defined for $t\neq s$. Sources (and sinks) are nevertheless very interesting for models of pedestrian, for instance in the case of people entering or exiting a door. It is interesting to recall that the $L^1$ distance (that one could try to use, since it is defined even between two measures with different masses) is not suitable in this context, since Lipschitzianity of $v$ with respect to $\mu$ measured in $L^1$ does not guarantee uniqueness (see \cite{nostro}).

In this article, to deal with a source in \r{e-cauchy}, we focus on the space of Borel measures with finite mass on $\R^d$ (denoted with $\Mu$), that we endow with the generalized Wasserstein distance $\gwbase$. We also denote with $\Muu$ the subspace of $\Mu$ of measures that are absolutely continuous with respect to the Lebesgue measure and with bounded support. In this framework, we prove existence and uniqueness of the solution of \r{e-cauchy} with $\mu_0\in\Muu$ under the following hypotheses:

\newcommand{\Hp}{{\bf (H)}}

\vspace{3mm}
\noindent\begin{Sbox}\begin{minipage}{\textwidth}\vspace{2mm}\begin{center}\Hp\vspace{5mm}\\
\begin{minipage}{0.9\textwidth}
The function
\mmfunz{v\Pq{\mu}}{\Mu}{C^{1}(\R^d)\cap L^\infty(\R^d)}{\mu}{v\Pq{\mu}}
satisfies 
\bi
\i $v\Pq{\mu}$ is uniformly Lipschitz and uniformly bounded, i.e. there exist $L$, $M$ not depending on $\mu$, such that for all $\mu\in\Mu, x,y\in\R^d,$
\bqn
\hspace{-5mm}|v\Pq{\mu}(x)-v\Pq{\mu}(y)|\leq L |x-y|\qquad |v\Pq{\mu}(x)|\leq M.
\eqnn

\i $v$ is a Lipschitz function, i.e. there exists $N$ such that 
\bqn\|v\Pq{\mu}-v\Pq{\nu}\|_{\mathrm{C^0}} \leq N \gw{\mu,\nu}.
\eqnn
\ei

The function
\mmfunz{h\Pq{\mu}}{\Mu}{\Muu}{\mu}{h\Pq{\mu}}
satisfies 
\bi
\i $h\Pq{\mu}$ has uniformly bounded mass and support, i.e. there exist $P,R$ such that
\bqn
 h\Pq{\mu}(\R^d)\leq P,\qquad \mathrm{supp} \Pt{h\Pq{\mu}}\subseteq B_R(0).
\eqnn

\i $h$ is a Lipschitz function, i.e. there exists $Q$ such that 
\bqn \gw{h\Pq{\mu},h\Pq{\nu}} \leq Q \gw{\mu,\nu}.
\eqnn
\ei
\end{minipage}\end{center}\vspace{2mm}\end{minipage}\end{Sbox}\fbox{\TheSbox}\vspace{3mm}

\brem
The hypotheses \Hp\ can be relaxed in standard ways. For example, the results still hold if we remove the uniform boundedness of $v$ on $\R^d$ and ask for uniform boundedness only in a point $x_0$. See e.g. \cite[Chap. 2]{libro-bressan}.
\erem

\brem The application to pedestrian dynamics also explains the choice of the basic assumptions \Hp, namely that we deal with measures with bounded support.
\erem

The structure of the paper is the following. In Section \ref{s-gw}, we define the generalized Wasserstein distance $\gwbase$ and we prove some important properties related to that. In particular, this distance metrizes the weak topology for tight sequences. Moreover, $\Mu$ is complete with respect to such distance. We also compare $\gwbase$ with other distances, like the Levy-Prokhorov distance. We then restrict ourselves to the study of $\gwbase$ in $\Muu$ and provide Gronwall-like estimates under flow action.

In Section \ref{s-equazione}, we describe the complete picture for \r{e-cauchy} under \Hp. We first provide a candidate solution for \r{e-cauchy} via a semi-discrete Lagrangian scheme and using the sample-and-hold method. We then show that it is indeed a solution, and finally that it is unique.

\section{Generalized Wasserstein distance}
\label{s-gw}

\newcommand{\Der}[2]{{D_{#2}#1}}

In this section, we define the generalized Wasserstein distance $\gw{\mu,\nu}$ that we study in this article, and prove some useful properties. We first recall basic definitions and notations about measure theory and Wasserstein distance. For a complete introduction, see \cite{ev-gar,villani}.

\subsection{Notation and standard Wasserstein distance}

In this section, we fix the notation that we use throughout the paper, and recall definitions and properties related to measure theory and the Wasserstein distance, like push-forward of measures $\gamma\#\mu$ and transference plans.

Let $\mu$ be a positive Borel measure with locally finite mass. If $\mu_1$ is absolutely continuous with respect to $\mu$, we write $\mu_1\ll \mu$. If $\mu_1\ll\mu$ and $\mu_1(A)\leq\mu(A)$ for all Borel sets, we write $\mu_1\sotto\mu$. Given a measure with finite mass, we denote with $|\mu|:=\mu(\R^d)$ its norm. More in general, if $\mu=\mu^+-\mu^-$ is a signed Borel measure, we have $|\mu|:=|\mu^+|+|\mu^-|$. Such norm defines a distance in $\Mu$, that is $|\mu-\nu|$. It is useful to recall that, if $\mu_1\ll\mu$ and $d\mu_1=f\, d\mu$ with $f\in L^1(d\mu)$, then $|\mu_1|=\int |f|\,d\mu$.

Given two measures $\mu,\nu$, one can always write in a unique way $\mu=\mu_{ac}+\mu_s$ such that $\mu_{ac}\ll \nu$ and $\mu_{s}\perp\nu$, i.e. there exists $B$ such that $\mu_s(B)=0$ and $\nu(\R^n\setminus B)=0$. This is the Lebesgue's decomposition Theorem. Then, it exists a unique $f\in L^1(d\nu)$ such that $d\mu_{ac}(x)=f(x)\,d\nu(x)$. Such function is called the Radon-Nikodym derivative of $\mu$ with respect to $\nu$. We denote it with $\Der{\mu}{\nu}$. For more details, see e.g. \cite{ev-gar}.\\

Given a Borel map $\fz{\gamma}{\R^d}{\R^d}$, one can consider the following action on a measure $\mu\in\Mu$:
\bqn
\gamma\#\mu(A):=\mu(\gamma^{-1}(A)).
\eqnn
An evident property is that the mass of $\mu$, i.e. $\mu(\R^d)$ is identical to the mass of $\gamma\#\mu$.

Given two measures $\mu,\nu$ with the same mass, it is thus possible to ask if there exists a $\gamma$ such that $\nu=\gamma\#\mu$. We say that $\gamma$ sends $\mu$ to $\nu$. Moreover, one can add a cost to such $\ga$, given by $I\Pq{\gamma}:=|\mu|^{-1}\,\int_{\R^d} |x-\ga(x)|^p \,d\mu(x)$. This means that each infinitesimal mass $\delta\mu$ is sent to $\delta \nu$ and that its infinitesimal cost is related to the $p$-th power of the distance between them. The, one can consider the map $\gamma$ minimizing such cost, if it exists. This is known as the Monge problem, stated by Monge in 1791.

In general, this procedure works only for special $\mu,\nu$ and $p$. Indeed, there exist simple examples of $\mu,\nu$ for which a $\gamma$ that sends $\mu$ to $\nu$ does not exist. For example $\mu=2\delta_1$, $\nu=\delta_{0}+\delta_2$ on the real line have the same mass, but there exists no $\gamma$ with $\nu=\gamma\#\mu$, since $\gamma$ cannot separate masses. Moreover, one can have a sequence $\gamma_n$ of maps such that $I\Pq{\gamma_n}$ is a minimizing sequence, but the limit is not a map $\gamma^*$. A simple condition that ensures the existence of a minimizing $\gamma$ is that $\mu$ and $\nu$ are absolutely continuous with respect to the Lebesgue measure.

For such reason, one can generalize the problem to the following setting. Given a probability measure $\pi$ on $\R^d\times\R^d$, one can interpret it as a method to transfer a measure $\mu$ on $\R^d$ to another measure on $\R^d$ as follows: each infinitesimal mass on a location $x$ is sent to a location $y$ with a probability given by $\pi(x,y)$. Formally, $\mu$ is sent to $\nu$ if the following properties hold:
\bqn
|\mu|\,\int_{\R^d} d\pi(x,\cdot)=d\mu(x),\qquad \qquad |\nu|\,\int_{\R^d} d\pi(\cdot,y)=d\nu(y).
\eqnl{e-pi}
Such $\pi$ is called a transference plan from $\mu$ to $\nu$. We denote the set of such transference plans as $\Pi(\mu,\nu)$. Since one usually deals with probability measures $\mu,\nu$, the terms $|\mu|,|\nu|$ are usually neglected in the literature. A condition equivalent to \r{e-pi} is that, for all $f,g\in C^\infty_c(\R^d)$ it holds $|\mu|\,\int_{\R^d\times \R^d} (f(x)+g(y))\,d\pi(x,y) = \int_{\R^d} f(x)\,d\mu(x)+ \int_{\R^d} g(y)\,d\nu(y)$. Then, one can define a cost for $\pi$ as follows $J\Pq{\pi}:=\int_{\R^d\times\R^d} |x-y|^p \,d\pi(x,y)$ and look for a minimizer of $J$ in $\Pi(\mu,\nu)$. Such problem is called the Monge-Kantorovich problem.

It is important to observe that such problem is a generalization of the Monge problem. Indeed, given a $\gamma$ sending $\mu$ to $\nu$, one can define a transference plan $\pi=(\Id\times \gamma)\#\mu$, i.e. $d\pi(x,y)=\mu(\R^n)^{-1}\, d\mu(x)\delta_{y=\gamma(x)}$. It also holds $J\Pq{\Id\times\gamma}=I\Pq{\gamma}$. The main advantages of such approach are the following: first, the existence of at least one $\pi$ satisfying \r{e-pi} is easy to check, since one can choose $\pi(A\times B)=|\mu|^{-1}\,\mu(A)\nu(B)$, i.e. the mass from $\mu$ is proportionally split to $\nu$. Moreover, a minimizer of $J$ in $\Pi(\mu,\nu)$ always exists.\\

A natural space in which $J$ is finite is the space of Borel measures with finite $p$-moment, that is
\bqn
\Mu^p:=\Pg{\mu\in\Mu\ |\ \int |x|^p\, d\mu(x)<\infty}. 
\eqnn
In the following, we also denote with $\P$ the space of probability measures, i.e. the measures in $\Mu$ with unit mass. We also deal with $\P^p:=\Mu^p\cap \P$, i.e. the space of probability measures with finite $p$-moment. One can thus define on $\Mu^p$ the following operator between measures of the same mass, called the \b{Wasserstein distance}:
\bqn
W_p(\mu,\nu)=(|\mu|\,\min_{\pi\in\Pi(\mu,\nu)} J\Pq{\pi})^{1/p}.
\eqnn
It is indeed a distance on the subspace of measures in $\Mu^p$ with a given mass, see \cite{villani}. It is easy to prove that $W_p(k\mu,k\nu)=k^{1/p}W_p(\mu,\nu)$ for $k\geq 0$, by observing that $\Pi(k\mu,k\nu)=\Pi(\mu,\nu)$ and that $J\Pq{\pi}$ does not depend on the mass. We will recall some other properties all along the paper, when useful.\\

\subsection{Definition of generalized Wasserstein distance}
We are now ready to define the \b{generalized Wasserstein distance} $\gw{\mu,\nu}$. We first give a rough description of the idea. Imagine to have three different admissible actions on $\mu,\nu$: either add/remove mass to $\mu$, or add/remove mass to $\nu$ or transport mass from $\mu$ to $\nu$. The three techniques have their cost: add/remove mass has a unitary cost $a$ (in both cases); transport of mass has the classic Monge-Kantorovich cost $J$, multiplied by a fixed constant $b$. The distance is the minimal cost of a mix of such techniques. We will show in the following that, depending on $\mu,\nu$, all mixes are possible: either remove all the mass of $\mu$ and $\nu$ (if they are very far), or transport the whole $\mu$ to $\nu$ (if they have the same mass and are close enough), or a mix of the two (for example when $\mu$ and $\nu$ are very close but with different masses). Instead, we will prove that add mass is never optimal.

\brem The fact that the unitary cost of adding/removing mass is identical for the two terms $\mu,\nu$ is to ensure symmetry of $\gwbase$.
\erem

We now formally define the generalized Wasserstein distance.
\bdeff
Let $\Mu$ be the space of Borel measures with finite mass on $\R^d$. Then, given $a,b\in(0,\infty)$ and $p\geq 1$, the \b{generalized Wasserstein distance} is
\bqn
\gw{\mu,\nu}=\inf_{\scriptstyle \ba{c}\tilde\mu,\tilde\nu\in\Mu^p\\
|\tilde\mu|=|\tilde\nu|\ea}\Pt{a|\mu-\tilde\mu|+a|\nu-\tilde\nu|+bW_p(\tilde\mu,\tilde\nu)}.
\eqnl{e-gw}
\edeff

\bp \label{p-distance}
The operator $\gwbase$ is a distance. Moreover, one can restrict the computation in \r{e-gw} to $\tilde\mu\sotto\mu,\tilde\nu\sotto\nu$ and the infimum is always attained.
\ep

To prove this Proposition, we need to recall some known concepts and results related to tightness and weak convergence. For more details, see \cite{old-new,villani}.

We first recall results about tightness. We always deal with measures on $\R^d$.
\bdeff \label{d-tight} A set of measures $M$ is tight if for each $\eps>0$ there exists a compact $K_\eps$ such that $\mu(\R^d\setminus K_\eps)<\eps$ for all $\mu\in M$.
\edeff
\bl The following holds: \be
\i A measure with finite mass is tight.
\i A set of measures $M$ such that all measures are bounded above by a tight measure $\tilde\mu$ is tight.
\i The union of a finite number of tight sets is tight. In particular, a finite set of measures with finite mass is tight.
\ee
\el
\bproof The first result is proved by taking a sequence of invading compacts $\Pg{K^n}$. Since $\mu(\R^d)=\lim_n \mu(K^n)=\sum_{n=1}^\infty \mu(K^n\setminus K^{n-1})<\infty$, then $\mu(\R^d\setminus K^i)= \sum_{n=i+1}^\infty \mu(K^n\setminus K^{n-1})$, as small as needed, since it is the remainder of a converging sequence.

For the second property, for each $\eps$ take $K_\eps$ that gives tightness of $\tilde \mu$. Since $\mu\in M$ implies $\mu\sotto\tilde\mu$, hence $\mu(\R^d\setminus K_\eps)<\eps$.

For the third result, take the sets $M_1,\ldots,M_n$, each of them being tight. For each $\eps>0$, there exists a corresponding compact $K^i_\eps$ such that for each $\mu\in M_i$ satisfies $\mu(\R^d\setminus K^i_\eps)$. Then, define $K_\eps=\cup_{i=1}^n K^i_\eps$, that is compact because it is a finite union of compacts. Now take $\mu\in \cup_{i=1}^n M_i$ and observe that $\mu\in M_i$ for some $i$, thus $\mu(\R^d\setminus K_\eps)\leq \mu(\R^d\setminus K_\eps^i)<\eps$.
\eproof

We now recall results about weak convergence and give some technical lemmas we use in the following.
\bt[Prokhorov's theorem] Let $X$ be a Polish space. A set $P$ in the space of probabilities $\mathcal{P}(X)$ is precompact for the weak topology if and only if it is tight.
\et
We recall that $\R^d$ is a Polish space (see a definition in \cite{villani}), thus Prokhorov's theorem can be applied in our setting.

\bt[Weak compactness \cite{ev-gar}] \label{t-weak} Let $\mu^n$ be a sequence of Radon measures in $\R^d$ with uniformly bounded mass\footnote{The result in \cite{ev-gar} is stated with a weaker condition, that is uniformly boundedness of mass on each compact.}, i.e. there exists $M$ such that $|\mu^k|\leq M$ for all $k$. Then there exists a subsequence $\mu^{n_j}$ and a Radon measure $\mu^*$ such that $\mu^{n_j}\weak \mu^*$.
\et
From now on, we denote with the asterisk $\mu^*,c^*$,... the limit (or the weak limit) of a sequence $\mu_n,c_n$,...

\bl \label{l-transf} Let $\mu^n,\nu^n$ be two sequences of Borel measures such that $\mu^n(\R^d)=\nu^n(\R^d)$ for each $n$. For each $n$, let $\pi^n$ be a transference plan with marginals $\mu^n,\nu^n$. If we have $\mu^n\weak \mu^*$, $\nu^n\weak \nu^*$, $\pi^n\weak \pi^*$ for some $\mu^*,\nu^*,\pi^*$, then $\pi^*$ is a transference plan with marginals $\mu^*,\nu^*$.
\el
\bproof The weak convergence of $\pi^n$ means that, for each $f\in C_c^{\infty}(\R^d\times \R^d)$ we have $\int f(x,y) \,d\pi^n(x,y) \to \int f(x,y) \,d \pi^*(x,y)$. In particular, choose $f(x,y)=g(x)$ and observe that $$\int g(x) \,d \pi^*(x,y)\leftarrow \int g(x) \,d\pi^n(x,y)=\int g(x) \,d\mu^n(x)\to \int g(x) \, d\mu^*(x).$$
By uniqueness of the limit, we have the proof for the marginal $\mu^*$. Using $f(x,y)=h(y)$, we have the same for $\nu^*$.
\eproof

We are now ready to prove Proposition \ref{p-distance}.

{{\it Proof of Proposition \ref{p-distance}.}} The \b{symmetry property} $\gw{\mu,\nu}=\gw{\nu,\mu}$ is evident.

We first prove that \b{we can always restrict to $\tilde\mu\sotto\mu$}. Define $$\mathcal{C}(\tilde\mu,\tilde\nu):=a|\mu-\tilde\mu|+a|\nu-\tilde\nu|+bW_p(\tilde\mu,\tilde\nu).$$ First assume that the infimum of $\mathcal{C}$ is attained by $\tilde\mu\not\sotto\mu$ and a certain $\tilde\nu$. Let $\pi\in\Pi(\tilde\mu,\tilde\nu)$ be the transference realizing $W_p(\tilde\mu,\tilde\nu)$. Let $d$ be the Radon-Nikodym derivative $f=\Der{\tilde\mu}{\mu}$ and $\mu_\perp:=\tilde\mu-f\mu$ the orthogonal of $\tilde\mu$ with respect to $\mu$. Define $\bar\mu:=\min\Pg{f,1}\mu$ and $\bar\nu$ the image of $\bar\mu$ under $\pi$. Since $\bar\mu\sotto \tilde\mu$ and $\pi\in\Pi(\tilde\mu,\tilde\nu)$, then $\bar\nu\sotto \tilde \nu$. Moreover, $|\tilde\nu-\bar\nu|=|\tilde\mu-\bar\mu|$, since $|\tilde\mu|=|\tilde\nu|$, $|\bar\mu|=|\bar\nu|$ by construction. Observe that 
\bqn
|\mu-\tilde\mu|&=&\int|1-f|\,d\mu+\mu_\perp(\R^d)=\int_{1-f\geq 0}(1-f)\,d\mu+\int_{1-f<0} (f-1)\,d\mu+\mu_\perp(\R^d)=\nn
&=&\int_{f\leq 1}(d\mu-d\bar\mu)+\Pt{\int_{f>1} (d\tilde\mu-d\bar\mu)-\mu_\perp(\Pg{f>1})}+\mu_\perp(\R^d)\geq|\mu-\bar\mu|+|\bar\mu-\tilde\mu|.
\eqnn
Hence
\bqn
|\mu-\bar\mu|+|\nu-\bar\nu|\leq |\mu-\tilde\mu|-|\bar\mu-\tilde\mu|+|\nu-\tilde\nu|+|\tilde\nu-\bar\nu|=|\mu-\tilde\mu|+|\nu-\tilde\nu|.
\eqnl{e-l1uguale}
The fact that $W_p(\bar\mu,\bar\nu)\leq W_p(\tilde\mu,\tilde\nu)$ is a direct consequence of the fact that $\pi\in \Pi(\tilde\mu,\tilde\nu)$ can be restricted to $\pi'\in \Pi(\bar\mu,\bar\nu)$, by construction of $\bar\nu$. Moreover, the cost of $\pi'$ is smaller\footnote{Theorem 4.6 in \cite{old-new} also shows that $\pi'$ is optimal, but this is not crucial here.} that the cost of $\pi$. Using this inequality and \r{e-l1uguale}, we have that $\mathcal{C}(\bar\mu,\bar\nu)\leq\mathcal{C}(\tilde\mu,\tilde\nu)$. Since $\bar\mu\sotto\mu$, then the result is proven.

If the infimum is not attained, consider a minimizing sequence $\tilde\mu^n$ and construct each $\bar\mu^n:=\min\Pg{\Der{\tilde\mu^n}{\mu},1}\mu$, that gives another minimizing sequence with $\bar\mu^n\sotto\mu$. By symmetry, the same property can be proved for the term $\nu$.

We now prove that \b{the infimum in \r{e-gw} is always attained}. To prove it, we restrict ourselves to $\tilde\mu\sotto\mu, \tilde\nu\sotto\nu$.   Take a sequence $\tilde\mu^n,\tilde\nu^n$ such that $\mathcal{C}(\tilde\mu^n,\tilde\nu^n)\to \gw{\mu,\nu}$ and $\tilde\mu^n\sotto \mu,\tilde\nu^n\sotto \nu$. Since both $\mu, \nu\in\Mu$ have finite mass, then $\Pg{\tilde\mu^n}$, $\Pg{\tilde\nu^n}$ have both uniformly bounded masses. Thus, due to Theorem \ref{t-weak}, passing to sub-sequences we have $\tilde\mu^n\weak \mu^*$, $\tilde\nu^n\weak  \nu^*$. We now prove that $\mathcal{C}(\mu^*,\nu^*)\leq\lim_n\mathcal{C}(\tilde\mu^n,\tilde\nu^n)=\gw{\mu,\nu}$. First recall that weak convergence gives $|\mu-\mu^*|\leq \liminf_n |\mu-\tilde\mu^n|$ and equivalently for $|\nu-\nu^*|$. We are left to prove that 
\bqn
W_p(\mu^*, \nu^*)\leq \lim_n W_p(\tilde\mu^n,\tilde\nu^n).
\eqnl{e-daprov}
If $\mu^*=\nu^*=0$, then we are done. Otherwise, the sequence $c_n:=|\tilde\mu^n|=|\tilde\nu^n|$ (eventually passing to a sub-sequence) converges to $c^*>0$. For $n$ such that $c_n\neq 0$, define the probability measures $\bar\mu^n=c_n^{-1} \tilde\mu^n$, $\bar\nu^n=c_n^{-1}\tilde\nu^n$. It is clear that $\bar\mu^n\weak \bar\mu^*=(c^*)^{-1} \mu^*$, and similarly for $\bar\nu^n$. Denote with $\pi^n$ the optimal transference plan in $\Pi(\tilde\mu^n,\tilde\nu^n)=\Pi(\bar\mu^n,\bar\nu^n)$. Since $\bar\mu^n\leq \frac{1}{\sup c_n} \mu$ and $\bar\nu^n\leq \frac{1}{\sup c_n} \nu$, then $M:=\Pg{\bar\mu^n}$, $N:=\Pg{\bar\nu^n}$ are both tight, hence the set of transference plans $\Pi(M,N)$ is tight (see e.g. \cite[Lemma 4.4]{old-new}). Hence, due to Prokhorov's theorem, up to sub-sequences we have $\pi^n\weak \pi^*$ for some $\pi^*$. Using Lemma \ref{l-transf}, we have that $\pi^*$ is  a transference plan with marginals $\bar\mu^*$ and $ \bar\nu^*$ (not necessarily optimal). Since the distance is non-negative, then the functional $J:\pi\to \int |x-y|^p\, d\pi(x,y)$ is lower semicontinuous with respect to the weak topology, see \cite[Lemma 4.3]{old-new}, thus 
\bqn
W_p(\mu^*,\nu^*)=(c^*)^{1/p} W_p(\bar\mu^*,\bar\nu^*)\leq (c^*)^{1/p} J(\pi^*)^{1/p} \leq \lim_n c_n^{1/p} J(\pi^n)^{1/p}= \lim_n c_n^{1/p} W_p(\bar\mu^n,\bar\nu^n)=\lim_n W_p(\tilde\mu^n,\tilde\nu^n).
\eqnn

We now prove that \b{$\gw{\mu,\nu}=0$ implies $\mu=\nu$}. Since the infimum is attained for some $\tilde\mu,\tilde\nu$, then $|\mu-\tilde\mu|=W_p(\tilde\mu,\tilde\nu)=|\nu-\tilde\nu|=0$ implies $\mu=\tilde\mu=\tilde\nu=\nu$.

We now prove \b{triangle inequality} $\gw{\mu,\eta}\leq \gw{\mu,\nu}+\gw{\nu,\eta}$.  Denote with $\tilde\mu,\tilde\nu^1$ the minimizers in \r{e-gw} for $\gw{\mu,\nu}$ and $\tilde\nu^2,\tilde\eta$ the minimizers for $\gw{\nu,\eta}$. Observe that we have $\tilde \nu^1\neq \tilde\nu^2$ in general. Also call $\pi^1, \pi^2$ the transference plans realizing $W_p(\tilde\mu,\tilde\nu^1),W_p(\tilde\nu^2,\tilde\eta)$, respectively. Define now $\bar \nu:=\min\Pg{\Der{\tilde\nu^1}{\tilde\nu^2},1}\tilde \nu_2$. Define $\bar\mu, \bar\eta$ the marginals of $\bar \nu$ with respect to $\pi^1,\pi^2$ respectively, i.e. $\pi^1\in \Pi(\bar\mu,\bar\nu)$ and $\pi^2\in\Pi(\bar\nu,\bar\eta)$. Thus we have $\pi^2\circ\pi^1\in\Pi(\bar\mu,\bar\eta)$. By construction, we have $|\tilde\mu-\bar\mu|=|\tilde\nu^1-\bar\nu|$. Moreover, $W_p(\bar\mu,\bar\nu)\leq W_p(\tilde\mu,\tilde\nu^1)$, because we use the same transference plan $\pi^1$ with a smaller mass. Similar properties hold for $\tilde\eta,\bar\eta$. We thus have
\bqn
\gw{\mu,\eta}&\leq& a|\mu-\bar\mu|+a|\eta-\bar\eta|+b W_p(\bar\mu,\bar\eta)\leq\nn
&\leq&a|\mu-\tilde\mu|+a|\tilde\mu-\bar\mu|+a|\eta-\tilde\eta|+a|\tilde\eta-\bar\eta|+b W_p(\bar\mu,\bar\nu)+bW_p(\bar\nu,\bar\eta)\leq\nn
&\leq&a|\mu-\tilde\mu|+a|\tilde\nu^1-\bar\nu|+a|\eta-\tilde\eta|+a|\tilde\nu^2-\bar\nu|+b W_p(\tilde\mu,\tilde\nu^1)+bW_p(\tilde\nu^2,\tilde\eta).
\eqnl{e-serve}

Define $f:=\Der{\tilde\nu^1}{\tilde\nu^2}$ and ${\tilde\nu^1}_\perp:=\tilde\nu^1-f\tilde\nu^2$ the orthogonal part of $\tilde\nu^1$ with respect to $\tilde\nu^2$. Observe that both $\tilde\nu^1,\tilde\nu^2\sotto \nu$, hence $\Der{\nu}{\tilde\nu^2}\geq \max\Pg{f,1}$ and $\nu_\perp:=\nu-(\Der{\nu}{\tilde\nu^2})\tilde\nu^2$ satisfies $\nu_\perp\geq {\tilde\nu^1}_\perp$. Then one has
\bqn
|\tilde\nu^1-\bar\nu|&=&\int \Pt{f-\min\Pg{f,1}} d\tilde\nu^2+{\tilde\nu^1}_\perp(\R^d)=\int_{f\geq 1} (f-1)d\tilde\nu^2+\int_{f< 1} 0 \,d\tilde\nu^2+{\tilde\nu^1}_\perp(\R^d)=\nn
&=&\int_{f\geq 1} (\max\Pg{f,1}-1)\,d\tilde\nu^2+{\tilde\nu^1}_\perp(\R^d)\leq\int_{f\geq 1} (d\nu-d\tilde\nu^2)+\nu_\perp(\R^d)= |\nu-\tilde\nu^2|,
\eqnn
and similarly $|\tilde\nu^2-\bar\nu|\leq |\tilde\nu-\tilde\nu^1|$. Plugging them into \r{e-serve}, one has the proof.
{\hfill$\square$}

One interesting feature of this distance is that the $|\cdot|$ term and the Wassertein term $W_p$ have different degree of homogeneity with respect to translation in $\R^n$, thus the optimal strategy for $\gwbase$ varies when translating one measure. For example, compute $\gw{\delta_0,\delta_x}$ as a function of $x\geq 0$. We have $|\delta_0-\delta_x|=2$ and $W_p(\delta_0,\delta_x)=x$. Hence $\gw{\delta_0,\delta_x}=\min\Pg{2a,bx}$. If $2a<bx$, i.e. measures are \virg{far}, then the optimal strategy is to delete both masses $\delta_0$ and $\delta_x$, otherwise it is optimal to move $\delta_0$ to $\delta_x$ with a translation.

\begin{figure}[htb]
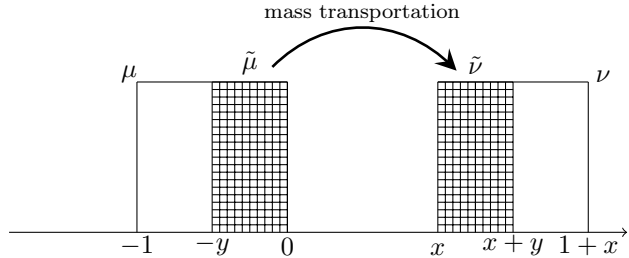


\begin{center}
\begin{pgfpicture}{0cm}{0cm}{9 cm}{4cm}

\begin{pgfscope}
\pgfsetendarrow{\pgfarrowto}
\pgfline{\pgfxy(0.3,0.3)}{\pgfxy(8.5,0.3)}
\end{pgfscope}

\pgfline{\pgfxy(2,0.3)}{\pgfxy(2,2.3)}
\pgfline{\pgfxy(4,2.3)}{\pgfxy(2,2.3)}
\pgfline{\pgfxy(4,2.3)}{\pgfxy(4,0.3)}
\pgfputat{\pgfxy(2,2.3)}{\pgfbox[right,bottom]{$\mu$}}
\pgfputat{\pgfxy(2,0)}{\pgfbox[center,bottom]{$-1$}}
\pgfputat{\pgfxy(4,0)}{\pgfbox[center,bottom]{$0$}}

\pgfline{\pgfxy(3,0.3)}{\pgfxy(3,2.3)}
\pgfgrid[stepx=1mm,stepy=1mm]{\pgfxy(3,0.3)}{\pgfxy(4,2.3)}
\pgfputat{\pgfxy(3.5,2.4)}{\pgfbox[center,bottom]{$\tilde\mu$}}
\pgfputat{\pgfxy(3,0)}{\pgfbox[center,bottom]{$-y$}}

\pgfline{\pgfxy(2+4,0.3)}{\pgfxy(2+4,2.3)}
\pgfline{\pgfxy(4+4,2.3)}{\pgfxy(2+4,2.3)}
\pgfline{\pgfxy(4+4,2.3)}{\pgfxy(4+4,0.3)}
\pgfputat{\pgfxy(8.1,2.3)}{\pgfbox[left,bottom]{$\nu$}}
\pgfputat{\pgfxy(6,0)}{\pgfbox[center,bottom]{$x$}}
\pgfputat{\pgfxy(8,0)}{\pgfbox[center,bottom]{$1+x$}}

\pgfline{\pgfxy(5+1.5,0.3)}{\pgfxy(5+1.5,2.3)}
\pgfgrid[stepx=1mm,stepy=1mm]{\pgfxy(6,0.3)}{\pgfxy(7,2.3)}
\pgfputat{\pgfxy(6.5,2.4)}{\pgfbox[center,bottom]{$\tilde\nu$}}
\pgfputat{\pgfxy(7,0)}{\pgfbox[center,bottom]{$x+y$}}

\begin{pgfscope}
\pgfsetlinewidth{1pt}
\pgfsetendarrow{\pgfarrowsingle}
\pgfxycurve(3.8,2.5)(4.5,3.2)(5.5,3.2)(6.2,2.5)
\pgfputat{\pgfxy(5,3.1)}{\pgfbox[center,bottom]{\footnotesize mass transportation}}
\end{pgfscope}

\end{pgfpicture}

\caption{Choice of $\tilde\mu,\tilde\nu$ (shaded) for the computation of $\gwbase$.}
\label{fig-esempio}
\end{center}

\end{figure}

Another simple example permits us to show that optimal strategies can be either based on removing mass only ($L^1$ strategy), or on transporting mass only ($W_p$ strategy), or by a mix of them. We give an example in Figure \ref{fig-esempio}. Take the measures $\mu,\nu$ on the real line with densities $d\mu=\Chi_{\Pq{-1,0}}\,d\lam$, $d\nu=\Chi_{\Pq{x,1+x}}\,d\lam$ with $x\geq 0$, where $\lam$ is the Lebesgue measure. It is clear that the optimal strategy amounts to choose $\tilde\mu$, $\tilde\nu$ with densities $\Chi_{\Pq{-y,0}}$, $\Chi_{\Pq{x,x+y}}$ respectively, for a certain parameter $y\in\Pq{0,1}$. The $L_1$ strategy is given by choosing $y=0$, while the $W_p$ strategy is given by $y=1$. We now prove that all values of $y$ can be optimal, depending on $x$ and the parameters $a,b,p$. We have $|\mu-\tilde\mu|=|\nu-\tilde\nu|=1-y$ and $W_p(\tilde\mu,\tilde\nu)=|y|^{1/p} (x+y)$. We choose for simplicity $a=b=p=1$. Thus $\gw{\mu,\nu}=\min_{y\in\Pq{0,1}}2-2y+xy+y^2$. A simple computation shows that the minimum is attained by $y=\frac{2-x}{2}$ if $x\in\Pq{0,2}$, and $y=0$ for $x\geq 2$. This clearly shows that if the measures are very close ($x=0$), then the best strategy is the $W_p$ (since $y=1$), while for measures that are far ($x$ big) the best strategy is the $L^1$. As stated above, varying $x\in\Pt{0,2}$ one has mixed strategies.

We now state some simple properties of $\gwbase$.
\bp
\label{p-semplici}
The following properties hold:
\bi
\i $\gw{k\mu,k\nu}\leq \max\Pg{k^{1/p},k}\gw{\mu,\nu}$ for $k\geq 0$,
\i $\gw{\mu_1+\mu_2,\nu_1+\nu_2}\leq \gw{\mu_1,\nu_1}+\gw{\mu_2,\nu_2}$.
\i $a \Pabs{|\mu|-|\nu|}\leq \gw{\mu,\nu}\leq a (|\mu|+|\nu|)$
\ei
\ep
\bproof
The first two properties are direct consequences of similar properties for $|\cdot|$ and $W_p$.

For the third, we first prove the inequality
$a \Pabs{|\mu|-|\nu|}\leq \gw{\mu,\nu}$. Without loss of generality, we assume $|\mu|\geq |\nu|$. Take any $\tilde\mu\sotto\mu$, $\tilde\nu\sotto \nu$ and observe that $|\mu-\tilde\mu|=|\mu|-|\tilde\mu|$, and similarly for $\nu,\tilde\nu$. Also recall that $|\tilde\mu|=|\tilde\nu|\leq |\nu|$ by construction. Now choose $\tilde\mu\sotto\mu,\tilde\nu\sotto\nu$ realizing $\gw{\mu,\nu}$ and observe that 
\bqn
\gw{\mu,\nu}\geq a|\mu-\tilde\mu|+a|\nu-\tilde\nu|=a\Pt{|\mu|-|\tilde\mu|+|\nu|-|\tilde\nu|}\geq a\Pt{|\mu|-|\nu|+0}.
\eqnn
We now prove the inequality $\gw{\mu,\nu}\leq a (|\mu|+|\nu|)$. Choose $\tilde\mu=\tilde\nu=0$ and observe that $\mathcal{C}(0,0)=a (|\mu|+|\nu|)$. Since $\gwbase$ is the infimum on all $\tilde\mu,\tilde\nu$, we have the inequality.
\eproof

\subsection{Topology of the generalized Wasserstein distance}

In this section, we prove that $\gwbase$ metrizes weak convergence for tight sequences. We also prove that $\R^d$ is complete with respect to $\gwbase$.

We first prove a simple lemma for $\gwbase$, stating that optimal choices $\tilde\mu,\tilde\nu$ are very close to each other in $\R^d$. The basic idea is that, if we want to transfer far mass between $\tilde\mu$ and $\tilde\nu$, then it is cheaper to remove such masses from both measures.

\bp \label{p-ottimi} Let $\mu,\nu\in\Mu$ and $\tilde\mu,\tilde\nu$ be the choices realizing $\gw{\mu,\nu}$. If $\tilde\mu$ has support contained in a compact $K$, then $\tilde\nu$ has bounded support, contained in the enlarged compact 
\bqn
K^d:=\cup_{x\in K} B(x,d)
\eqnl{e-enlarged}
with $d=\frac{2a}{b}$. Here we denote with $B(x,r)$ the closed ball centered in $x$ with radius $r$.

Similarly, let $\tilde\mu,\tilde\nu$ be the choices realizing $\gw{\mu,\nu}$, and $\pi$ the transference plan realizing $W_p(\tilde\mu,\tilde\nu)$. Let $\mu'\sotto \tilde\mu$ have  support contained in a compact $K$. Then the corresponding marginal $\nu'$ with respect to $\pi$ has compact support contained in $K^d$.
\ep
\bproof
We prove the second statement, since the first can be recovered by choosing $\mu'=\tilde\mu$. 

With the notations given in the statement, we prove that the support of $\nu'$ is contained in $K^d$, by contradiction. Assume that there exists $d'>d$ and a Borel set $D\subset \R^d\setminus K^{d'}$ such that $\nu'(D)>0$. Thus define $\bar\nu:=k\nu'_{|_{D}}$ with $k\leq 1$ chosen in such a way that $\bar\nu(D)<1$. Observe that $\bar\nu\sotto\nu'\sotto\tilde \nu$. Define $\nu^*:=\tilde\nu-\bar\nu\leq \tilde\nu$ and $\mu^*$ the corresponding marginal given by $\pi$. We now prove that the choice $\mu^*,\nu^*$ in \r{e-gw} gives a cost that is strictly less than the optimal choice $\tilde\mu,\tilde\nu$. First observe that
\bqn
 |\mu-\mu^*|+ |\nu-\nu^*| &=& |\mu-\tilde\mu|+|\tilde\mu-\mu^*|+ |\nu-\tilde\nu| +|\tilde\nu-\nu^*|=\nn
 &=&|\mu-\tilde\mu|+ |\nu-\tilde\nu| +2|\tilde\nu-\nu^*|=|\mu-\tilde\mu|+ |\nu-\tilde\nu| +2|\bar\nu|.
\eqnl{e-una}
Now observe that $W_p^p(\tilde\mu,\tilde\nu)=W_p^p(\mu^*,\nu^*)+W_p^p(\tilde\mu-\mu^*,\tilde\nu-\nu^*)$, by construction of $\mu^*,\nu^*$ via the optimal transference plan $\pi$. Observe that $\supp(\tilde\mu-\mu^*)\subseteq \supp(\mu')\subseteq K$ and that $\supp(\tilde\nu-\nu^*)=\supp(\bar\nu)\subseteq D\subseteq \R^d\setminus K^{d'}$. In particular, if $x\in \supp(\tilde\mu-\mu^*),y\in \supp(\tilde\nu-\nu^*)$, then $|x-y|\geq d'$. Thus, given the optimal transference plan $\pi'\in \Pi(\tilde\mu-\mu^*,\tilde\nu-\nu^*)$, we have
\bqn
W_p^p(\tilde\mu-\mu^*,\tilde\nu-\nu^*)= \int |x-y|^p d\pi'(x,y)\geq\int d'^p \,d(\tilde\nu-\nu^*)=|\bar\nu| d'^p.
\eqnl{e-dua}
Putting together \r{e-una} and \r{e-dua}, and using $|\bar \nu|^{1/p}\geq |\bar\nu|$ since $|\bar\nu|\leq 1$, we have
\bqn
\mathcal{C}\Pt{\mu^*,\nu^*}&\leq& a |\mu-\tilde\mu|+ a|\nu-\tilde\nu| +2a|\bar\nu| + 
bW_p(\tilde\mu,\tilde\nu)-bW_p(\tilde\mu-\mu^*,\tilde\nu-\nu^*)\leq\nn
 &\leq& a |\mu-\tilde\mu|+ a|\nu-\tilde\nu| + 
bW_p(\tilde\mu,\tilde\nu) +2a|\bar\nu|-bd' |\bar\nu|^{1/p}\leq \mathcal{C}(\tilde\mu,\tilde\nu)-b \Pt{d'-2\frac{a}{b}}|\bar\nu|<\mathcal{C}(\tilde\mu,\tilde\nu).
\eqnn
Thus $\tilde\mu,\tilde\nu$ is not optimal. Contradiction. Thus $\supp(\nu')\subseteq K^{d'}$ for all $d'>d$, hence $\supp(\nu')\subseteq K^d$.
\eproof

We now prove the following convergence theorem, stating that $\gwbase$ metrizes weak convergence for tight sequences.
\bt \label{t-convergence}
Let $\Pg{\mu_n}$ be a sequence of measures in $\R^d$, and $\mu_n,\mu\in\Mu$. Then
$$\gw{\mu_n,\mu} \to 0\mbox{~~~~~~is equivalent to~~~~~~}\mu_n\weak \mu \mbox{~~and~~}\Pg{\mu_n}\,\mbox{is tight}.$$
\et
\bproof
Fix the following notation: for each $n$, let $\tilde\mu_n,\tilde\nu_n$ be optimal choices in \r{e-gw} for $\gw{\mu_n,\mu}$ with $\tilde\mu_n\sotto \mu_n,\tilde\nu_n\sotto\mu$, and $\pi_n$ the transference plan realizing $W_p(\tilde\mu_n,\tilde\nu_n)$.

We first prove $\Leftarrow$. Fix $\eps>0$. Let $N$ be such that $\gw{\mu_n,\mu}<\eps$ for all $n>N$. Since $\Pg{\mu_n}$ is tight, then the set $M:=\Pg{\mu_n}\cup \Pg{\mu}$ is tight too. Given $\delta>0$, consider the corresponding $K_\delta$ giving tightness of $M$. By definition, all $m\in M$ satisfy $m(\R^d\setminus K_\delta)<\delta$. Due to weak convergence of $\mu_n$ to $\mu$, we also have $\mu_n(K_\delta)\to\mu(K_\delta)$. Choose $N'$ such that $\Pabs{|\mu_n(K_\delta)|-|\mu(K_\delta)|}<\eps$ for all $n>N'$. It means that there exists $\nu_n$ (positive or negative, supported in $K_\delta$) such that $(\mu_n+\nu_n)(K_\delta)=\mu(K_\delta)$. Define $\tilde\mu_n:=(\mu_n+\nu_n)_{|_{K_\delta}}$ and $\tilde\mu:=\mu_{|_{K_\delta}}$. It is clear that $\nu_n\weak 0$, hence $\tilde\mu_n\weak \tilde\mu$. Since $K_\delta$ has bounded diameter  and $W_p$ metrizes weak convergence in bounded spaces (see \cite[7.12]{villani}, recalled below in Theorem \ref{t-vil}), then $W_p(\tilde\mu_n,\tilde\mu)\to 0$. Take $N>N'$ such that $W_p(\tilde\mu_n,\tilde\mu)<\delta$ and $|\nu_n|<\delta$ for all $n>N$. We now estimate 
\bqn
\gw{\mu_n,\mu}&\leq&\mathcal{C}(\tilde\mu_n,\tilde\mu)=a|\mu_n-\tilde\mu_n|+a|\mu-\tilde\mu|+bW_p(\tilde\mu_n,\tilde\mu)<\nn
&<& a|\mu_n(\R^d\setminus K_\delta)|+a|\nu_n|+a|\mu(\R^d\setminus K_\delta)|+b\delta<(3a+b)\delta.
\eqnn
Choose $\delta=\eps/(3a+b)$ and have the result.

We now prove $\Rightarrow$. As a first, main step, we prove that \b{$\gw{\mu_n,\mu}\to 0$ implies that $\Pg{\mu_n}$ is tight}. Fix now $\eps>0$. There exists a corresponding $N$ such that $\gw{\mu_n,\mu}<\eps$ for all $n>N$, and a compact $K_\eps$ such that $\mu(\R^d\setminus K_\eps)<\eps$, since $\mu$ is tight. Define $\nu'_n:=\tilde\nu_n$ restricted to $K_\eps$, and $\mu'_n$ the corresponding measure given by the transference plan $\pi_n$ realizing $W_p(\tilde\mu_n,\tilde\nu_n)$. Using Proposition \ref{p-ottimi}, we have that the support of $\mu'_n$ is contained in $K_\eps^d$ with $d=2\frac{a}{b}$. Observe that such set does not depend on $n$. By construction, we have
\bqn
|\tilde\mu_n-\mu'_n|=|\tilde\mu_n|-|\mu'_n|=|\tilde\nu_n|-|\nu'_n|=\tilde\nu_n(\R^d\setminus K_\eps)<\eps.
\eqnn We now estimate
\bqn
\mu_n(\R^d\setminus K^d_\eps)&\leq& |\mu_n-\tilde\mu_n|+\tilde\mu_n(\R^d)-\tilde\mu_n(K^d_\eps)\leq \frac{1}{a}\gw{\mu_n,\mu}+|\tilde\mu_n|-|\mu'_n|\leq\frac{\eps}{a} +\eps.
\eqnn
This is the tightness of $\Pg{\mu_n}_{n>N}$. Observe now that the finite set of measures $\Pg{\mu_1, \ldots, \mu_N}$ with finite masses is tight. Since a finite union of tight sets is tight, we have that $\Pg{\mu_n}$ is tight.

We now observe that $|\mu_n|\leq |\mu|+\frac{1}{a}W_p(\mu_n,\mu)$ by Proposition \ref{p-semplici}, that is a converging sequence of real numbers, thus $\mu_n$ have uniformly bounded mass. We now apply Theorem \ref{t-weak}, that gives \b{weak convergence} of $\mu_n$ to a certain $\mu^*$. We prove that $\mu=\mu^*$. Using the implication $\Leftarrow$, we have $\gw{\mu_n,\mu^*}\to 0$. Thus $\gw{\mu,\mu^*}\leq \lim_n \gw{\mu,\mu_n}+\gw{\mu_n,\mu^*}=0$, hence $\mu=\mu^*$. \eproof

It is interesting to compare such result with a similar result for the standard Wasserstein distance, that we recall here.
\bt{\cite[7.12]{villani}}\label{t-vil} Let ${\mu_k}$ be a sequence of probability measures in $\P^p$, and $\mu\in\P$. Then the following statements are equivalent:
\bi
\i $W_p(\mu_k,\mu)\to 0$;
\i $\mu_k\weak \mu$ and the following condition holds
\bqn
\lim_{R\to\infty} \limsup_{k} \int_{|x|>R} |x|^p\, d\mu_k(x)=0;
\eqnl{e-vil}
\i $\mu_k\weak \mu$ and the $p$-moment converges, i.e.
\bqn \int |x|^p\, d\mu_k(x) \to \int |x|^p\, d\mu(x).
\eqnl{e-vil2}
\ei
\et

Condition \r{e-vil} is called \virg{tightness} condition, a notation that could create some confusion with respect to Definition \ref{d-tight}. Anyway, we remark that condition \r{e-vil} is stronger than Definition \ref{d-tight}, thus our Theorem \ref{t-convergence} applies on a wider class than Theorem \ref{t-vil}. First of all, Theorem \ref{t-vil} applies for measures that have all the same mass, otherwise $W_p(\mu_k,\mu)$ is not defined. Moreover, even taking a sequence of probability measures $\mu_k$, one can have convergence in $\gwbase$ and no convergence in $W_p$. This occurs exactly in the case in which $\Pg{\mu_k}$ is tight according to Definition \ref{d-tight} and not according to condition \r{e-vil}. For example, take the following sequence of probability measures: $\mu_k:=(1-k^{-p})\delta_0+k^{-p}\,\delta_k$. It is clear that such sequence converges weakly to $\mu^*=\delta_0$, and that $W_p(\mu_k,\mu^*)=W_p(k^{-p}\, \delta_0,k^{-p}\,\delta_k)=k^{-1} \int k^p\, \delta_0 =k^{p-1}\not\to 0$. Indeed, condition \r{e-vil2} is not satisfied, since the p-moment of $\mu^*$ is 0, while the $p$-moment of $\mu_k$ is $1$. Similarly, condition \r{e-vil} is not satisfied, since for each $R$ the measures $\mu_{k}$ with $k>R$ satisfy $\int_{|x|>R} d\mu_k=k^{-p}\int_{|x|>R} k^p \,\delta_k=1\not\to 0$. Instead, estimate $\gw{\mu_k,\mu^*}$ by choosing $\tilde\mu_k=\mu_k-k^{-p}\,\delta_k$, $\tilde\mu^*=\mu^*-k^{-p}\,\delta_0$, that gives $\gw{\mu_k,\mu^*}=a|k^{-p}\,\delta_k|+a|k^{-p}\,\delta_0|=2k^{-p}\to 0$. Thus, $\Pg{\mu_k}$ is tight according to Definition \ref{d-tight}, as it is easy to prove by observing that $\mu_k(\R\setminus \Pq{-n,n})< n^{-p}$.

It is clear that, on the contrary, if $\Pg{\mu_k}$ satisfies \r{e-vil}, then it satisfies Definition \ref{d-tight}. Indeed, let $\Pg{\mu_k}$ satisfy \r{e-vil}. For each $\eps>0$ there exist $R,m$ such that $\int_{|x|>R} |x|^p \,d\mu_k<\eps$ for all $k>m$. One can always assume $R>1$. Thus $\mu_k(\R^n\setminus B(0,R))<\eps$ for all $k>m$, that ensures tightness of $\Pg{\mu_k}$.\\

We now prove completeness of $\Mu$.

\bp $\Mu$ is complete with respect to $\gwbase$.
\ep
\bproof
Take a sequence $\mu_n$ that is Cauchy with respect to $\gwbase$. We first show that $\Pg{\mu_n}$ is tight. Fix $\delta>0$ and let $N$ be such that $\gw{\mu_n,\mu_{n+k}}<\delta$ for all $n\geq N, k\geq 0$. Fix the notation $\tilde\mu_k,\tilde\nu_k$ to denote the choices in \r{e-gw} realizing the minimizer $\gw{\mu_N,\mu_{N+k}}$, and $\pi_k$ the transference plan realizing $W_p(\tilde\mu_k,\tilde\nu_k)$. Let $K_\delta$ be such that $\mu_{N}(\R^d\setminus K_\delta)<\delta$, that exists since $\mu_{N}$ has finite mass. Since $\tilde\mu_k\sotto\mu_N$, then $\tilde\mu_k(\R^d\setminus K_\delta)<\delta$. Define now $\mu'_k$ the restriction of $\tilde\mu_k$ to $K_\delta$ and $\nu'_k$ the corresponding marginal with respect to $\pi_k$. Using Proposition \ref{p-ottimi}, we have that $\supp(\nu'_k)\subseteq K_\delta^d$ with $d=\frac{2a}{b}$.

As a consequence $\tilde\nu_k(K_\delta^d)\geq \nu'_k(K_\delta^d)=\mu'_k(K_\delta)=\tilde\mu_k(K_\delta)$. Since $\tilde\mu_k(\R^d)=\tilde\nu_k(\R^d)$ by construction, we have $\tilde\nu_k(\R^d\setminus K_\delta^d)\leq \tilde\mu_k(\R^d\setminus K_\delta)<\delta$. Since $a|\mu_{N+k}-\tilde\nu_k|\leq\gw{\mu_{N},\mu_{N+k}}<\delta$, we have that $\mu_{N+k}(\R^d\setminus K_\delta^d)\leq \tilde\nu_k(\R^d\setminus K_\delta^d) +\delta/a<\delta+\delta/a$. Choose $\delta =\frac{\eps}{1+1/a}$ and $K'_\eps=K_\delta^d$, then $\mu_n(\R^d\setminus K'_\eps)<\eps$ for $n\geq N$. Thus, we have tightness for $n\geq N$. Since $\Pg{\mu_1,\ldots,\mu_{N}}$ is a finite family of measures with finite mass, then it is tight, hence the whole $\Pg{\mu_n}$ is tight.

Observe that $\Pg{\mu_n}$ has also uniformly bounded mass, thus there exists a subsequence $\mu_{n_k} \weak \mu^*$ for a certain $\mu^*$, due to Theorem \ref{t-weak}. Using Theorem \ref{t-convergence}, we have that $\gw{\mu_{n_k},\mu^*}\to 0$, and by triangular inequality we have $\gw{\mu_n,\mu^*}\to 0$.
\eproof

\subsection{Comparison with other distances}

In this section, we compare $\gwbase$ with two distances proposed in the literature, namely the Levy-Prokhorov distance (see e.g. \cite{villani}) and the distance $W_{b_2}$ defined in \cite{fg} by A. Figalli and N. Gigli.

\newcommand{\dlp}{d_{LP}}
We first recall the definition of the Levy-Prokhorov distance $\dlp$ between two probability measures $\mu,\nu$:
\bqn
\dlp(\mu,\nu):=\inf\Pg{\al>0 \mbox{~~such that~for any closed $A$ it holds~~} \mu(A)\leq \nu(A^\al)+\al},
\eqnn
where $A^\al$ is the enlarged set $A^\al:=\cup_{x\in A} B(x,\al)$. It is clear that the two terms $\nu(A^\al)$ and $\al$ represent a $W_p$ perturbation and a $L_1$ perturbation, respectively. For such reason, there are some common ideas between $\dlp$ and $\gwbase$.

The main difference here is that $\dlp$ was defined for probability measures only, while we deal with measures with different (finite) masses. Nevertheless, even restricting to the space of probability measures, the two distances have different values. We study a remarkable case, that is the distance between $\mu=\delta_0$ and $\nu=\frac{1}{2} \delta_{-d_1}+\frac{1}{2} \delta_{d_2}$ on the real line. We study the values of $\dlp$ and $\gwbase$ as functions of $d_1,d_2$. For both distances, the goal is to choose the optimal coupling, on one side between $\frac{1}{2}\delta_0$ and $\frac{1}{2} \delta_{-d_1}$, and on the other side between $\frac{1}{2}\delta_0$ and $\frac{1}{2} \delta_{d_2}$.

We now compare $\dlp$ with $\gwbase$ of parameters $a=\frac 12, b=1$. Without loss of generality, we assume $d_1\leq d_2$. The set of different cases is the following:
\bi
\i The two masses $\frac{1}{2} \delta_{-d_1},\frac{1}{2} \delta_{d_2}$ are both \virg{far} from $\delta_0$. The best choice for both couplings  is to focus on $L_1$ distance. For the $\dlp$, it means that $d_1\geq 1$. In this case $\dlp(\mu,\nu)=1$. For $\gwbase$, it means that $d_1\geq \frac{2a}{b}=1$.  In this case $\gw{\mu,\nu}=\mathcal{C}\Pt{0,0}=1$.
\i One of the masses is \virg{close} and the other is \virg{far}. This means that $d_1\leq 1\leq d_2$. In this case, we use the $W_p$ distance for the coupling $\frac{1}{2}\delta_0,\frac{1}{2} \delta_{-d_1}$ and the $L_1$ distance for the coupling $\frac{1}{2}\delta_0,\frac{1}{2} \delta_{d_2}$. Thus $\dlp(\mu,\nu)=\sup\Pg{\frac12,d_1}$ and $\gw{\mu,\nu}=\mathcal{C}\Pt{\frac{1}{2}\delta_0,\frac{1}{2} \delta_{-d_1}}=\frac12+2^{-1/p}d_1$.
\i Both masses are \virg{close, but one is not very close}. This is true for $\frac{1}{2}\leq d_2\leq 1$. In this case, the optimal strategies are different for the two distances:
\bi
\i for $\dlp$, we use the $W_p$ distance for the coupling $\frac{1}{2}\delta_0,\frac{1}{2} \delta_{-d_1}$ and the $L_1$ distance for the coupling $\frac{1}{2}\delta_0,\frac{1}{2} \delta_{d_2}$. Thus $\dlp(\mu,\nu)=\sup\Pg{\frac12,d_1}$.
\i for $\gwbase$ we use the $W_p$ distance for both couplings. This gives $\gw{\mu,\nu}=C(\mu,\nu)=W_p(\mu,\nu)=\Pt{\frac{d_1^p+d_2^p}2}^{1/p}$.
\ei
\i Both masses are \virg{very close}. This is the case of $d_2\leq \frac12$. For both distance the best choice is to use the $W_p$ distance. Thus $\dlp(\mu,\nu)=d_2$ and $\gw{\mu,\nu}=C(\mu,\nu)=W_p(\mu,\nu)=\Pt{\frac{d_1^p+d_2^p}2}^{1/p}$.
\ei
It is easy to prove that the two distances $\dlp,\gwbase$ are equivalent as norms on the space $(d_1,d_2)$. Moreover, they have the same level lines around $d_1=0,d_2=0$ when choosing $p=\infty$.\\

We now recall the distance $W_{b_2}$ defined in \cite{fg}. Such distance is defined on a subset $\Omega$ with non-empty boundary $\partial\Omega$. The idea is that such boundary is an infinite reserve of mass, in the following sense. Given two measures $\mu,\nu$, even with different masses, one can either send a part of mass of $\mu$ to $\nu$ or to the boundary $\partial \Omega$. Similarly, the mass of $\nu$ that does not receive mass from $\mu$ goes to the boundary $\partial \Omega$. In both cases, the cost is computed via the Wasserstein distance, either from $\mu$ to $\nu$ or from $\mu$ to $\partial \Omega$.

This cost is based on an approach that is rather different than ours. First of all, it deals with a space $\Omega$ that has boundary, otherwise one could not deal with measures with different masses. Moreover, the cost of sending mass to the boundary (that is similar to delete mass in our approach) is computed as a Wasserstein distance and not with $L_1$, like in our case.

\subsection{Estimates of generalized Wasserstein distance under flow actions}

In this section, we study properties of $\gwbase$ when restricted to $\Muu$. In particular, we are interested in estimates about the variation of $\gw{\mu,\nu}$ under action of flows on $\mu,\nu$. These Gronwall-like properties will be useful for the study of solutions of \r{e-cauchy}.

We first recall a connection between flows actions on measures and transport equation. Take a Lipschitz vector field $v$, that generates a flow $\Phi^v_t$ for $t\geq 0$. The flow is a diffeomorphism for the space $\R^d$, thus one can define $\mu_t:=\Phi^v_t\#\mu_0$ for a given measure $\mu_0$. One has the following theorem.

\bt \label{t-misuretrasporto} Take a Lipschitz vector field $v$, and the flow $\Phi^v_t$ it generates. Given $\mu_0\in\Muu$, and $\mu_t:=\Phi^v_t\#\mu_0$, then $\mu=\mu_{\Pq{0,T}}$ is the unique solution of the linear transport equation
\bqn
\begin{cases}
\partial_t\mu_t+\nabla\cdot(v \mu_t)=0\\
\mu_{|_{t=0}}=\mu_0
\end{cases}
\label{e-misuretrasporto}
\eqn
in $C\Pt{[0,T],\Muu}$, where $\Muu$ is endowed with the weak topology.
\et
\bproof The proof is a direct consequence of {\cite[Thm 5.34]{villani}}. See details in \cite{nostro}.\eproof


To study \r{e-cauchy} in the setting of the generalized Wasserstein distance, it is interesting to check the evolution of $\gwbase$ under flow action. We first recall here some properties about standard Wasserstein distance, that we proved in \cite{nostro}. Here, we emphasize the fact that the Wasserstein distance is computed between two measures $\mu,\nu$ with same mass, in general different than 1.
\bp
\label{p-old}
Given $v,w$ bounded and Lipschitz vector fields of Lipschitz constant $L$ and $\mu,\nu\in\Muu$, the following holds:
\be 
\i $W_p(\Phi^v_t\#\mu,\Phi^v_t\#\nu)\leq e^{\frac{p+1}{p}Lt}W_p(\mu,\nu)$,
\i $W_p(\mu,\Phi^v_t\#\mu)\leq t\,\|v\|_{C^0}  \mu(\R^d)^{1/p}$,
\i $W_p(\Phi^v_t\#\mu,\Phi^w_t\#\nu)\leq e^{\frac{p+1}{p}Lt}W_p(\mu,\nu)+\mu(\R^d)^{1/p} \frac{e^{L t/p}(e^{Lt}-1)}L \|v-w\|_{C^0}$.
\ee
\ep

We now prove similar properties for the generalized Wasserstein distance.
\bp
\label{p-utile}
Given $v,w$ bounded and Lipschitz vector fields of Lipschitz constant $L$, the following holds:
\be 
\i $\gw{\Phi^v_t\#\mu,\Phi^v_t\#\nu}\leq e^{\frac{p+1}{p}Lt}\gw{\mu,\nu}$,
\i $\gw{\mu,\Phi^v_t\#\mu}\leq t\, \|v\|_{C^0} \mu(\R^d)^{1/p}$,
\i $\gw{\Phi^v_t\#\mu,\Phi^w_t\#\nu}\leq e^{\frac{p+1}{p}Lt}\gw{\mu,\nu}+\mu(\R^d)^{1/p} \frac{e^{L t/p}(e^{Lt}-1)}L \|v-w\|_{C^0}$.
\ee
\ep
\bproof
For the first property, take $\tilde\mu,\tilde\nu$ realizing \r{e-gw} for $\gw{\mu,\nu}$ with $\tilde\mu\sotto\mu,\tilde\nu\sotto\nu$. Then
\bqn
\gw{\Phi^v_t\#\mu,\Phi^v_t\#\nu}&\leq& a|\Phi^v_t\#\mu-\Phi^v_t\#\tilde\mu|+a|\Phi^v_t\#\nu-\Phi^v_t\#\tilde\nu|+bW_p(\Phi^v_t\#\tilde\mu,\Phi^v_t\#\tilde\nu)\leq \nn
&\leq& a|\mu-\tilde\mu|+a|\nu-\tilde\nu|+b e^{\frac{p+1}{p}Lt} W_p(\tilde\mu,\tilde\nu)\leq
e^{\frac{p+1}{p}Lt} \Pt{ a|\mu-\tilde\mu|+a|\nu-\tilde\nu|+b  W_p(\tilde\mu,\tilde\nu)}.
\eqnn
The proofs of the second and the third properties are equivalent, based on proofs of Proposition \ref{p-old} given in \cite{nostro}.
\eproof

\section{Existence and uniqueness of solutions of \r{e-cauchy}}
\label{s-equazione}

In this section we prove the existence and uniqueness of the solution of \r{e-cauchy}, under the hypothesis \Hp. The key tool is the construction of a candidate solution by sample-and-hold; we then prove that it is indeed a solution, and finally prove that it is unique.

\subsection{Construction of an approximated solution}

In this section, we apply the sample-and-hold method to construct a sequence of functions in $C\Pt{[0,T],\Muu}$ such that the limit exists and it is a solution of \r{e-cauchy}. With no loss of generality, we assume that the $T=1$.

We first define an approximated solution $\mu^k$ for each $k\in\N$. Given a fixed $k$, define $\dt:=\frac{1}{2^{k}}$ and consider the decomposition of the time interval in $\Pq{0,\dt},\Pq{\dt,2\dt},\Pq{2\dt,3\dt},\ldots,\Pq{(2^k-1)\dt,2^k\dt}$.

The idea is based on a Lagrangian scheme with the time-discretization provided above. In particular, we define a solution for each time $n\dt$, then consider that $v$ and $h$ are computed at this time and fixed on the next interval $[n\dt,(n+1)\dt)$. More precisely, for a fixed $k$, we define

\bi
\i $\mu_0^k:=\mu_0$;
\i $\mu_{(n+1)\dt}^k:=\Phi^{v\Pq{\mu^k_{n\dt}}}_\dt\#\mu^k_{n\dt} + \dt\,h\Pq{\mu^k_{n\dt}}$;
\i $\mu_t^k:=\Phi^{v\Pq{\mu^k_{n\dt}}}_\tau\#\mu^k_{n\dt} + \tau\,h\Pq{\mu^k_{n\dt}}$ with $n$ the maximum integer such that $t-n\dt\geq 0$ and $\tau:=t-n\dt$.
\ei

We now prove that, given $v,h$ satisfying \Hp, the sequence $\mu^k$ is a Cauchy sequence for the space $(C(\Pq{0,1},\Mu),\mathcal{D})$, where $$\mathcal{D}(\mu,\nu):=\sup_{t\in\Pq{0,1}}\gw{\mu_t,\nu_t}.$$ Since we have proved that $(\Mu,\gwbase)$ is a complete space, then also $(C(\Pq{0,1},\Mu),\mathcal{D})$ is complete. Observe that we deal with the whole $\Mu$ to have completeness, and that we will subsequently prove that the limit is indeed an element of $\Muu$.

\newcommand{\PP}{m}
\newcommand{\dtd}{{\scriptstyle \frac{\dt}{2}}}
We first make three simple observation:
\bi
\i At each time $t\in\Pq{0,1}$, the mass $\mu^k_t(\R^d)$ is bounded. Indeed, first observe that 
\bqn
\mu^k_{(n+1)\dt}(\R^d)& \leq & \Phi^{v\Pq{\mu^k_{n\dt}}}_\dt\#\mu^k_{n\dt}(\R^d) + \dt\,h\Pq{\mu^k_{n\dt}}(\R^d)\leq \mu^k_{n\dt}(\R^d)+\dt P.
\eqnn
As a consequence, we have $$\mu^k_t(\R^d)^{1/p}\leq \PP,\mbox{~~ with~~}m:=(\mu_0(\R^d)+ P)^{1/p}.$$ This holds for all $k,t$. We use such constant in the following estimates.
\i Given a fixed $k$, the evolution of the scheme satisfies $\gw{\mu_{(n+1)\dt}^k,\mu^k_{n\dt}}\leq \dt M m+\dt P$, and more in general
\bqn
\gw{\mu_{t}^k,\mu^k_s}\leq|t-s|(Mm+P).
\eqnn
\i We will use in the following the estimate
\bqn
\gw{\mu^{k+1}_{(n+\frac12)\dt},\mu^k_{n\dt}}&\leq& \gw{\mu^{k+1}_{(n+\frac12)\dt},\mu^{k+1}_{n\dt}}+\gw{\mu^{k+1}_{n\dt},\mu^k_{n\dt}}\leq\nn
&\leq& \dtd(Mm+P)+\gw{\mu^{k+1}_{n\dt},\mu^k_{n\dt}},
\eqnl{e-rough}
that is a simple consequence of Proposition \ref{p-utile}-2.
\ei

We now estimate the distance $\gw{\mu^{k+1}_{(n+1)\dt},\mu^k_{(n+1)\dt}}$ with respect to $\gw{\mu^{k+1}_{n\dt},\mu^k_{n\dt}}$, i.e. the evolution of $\gwbase$ at the discretization points of the Lagrangian scheme for $\mu^k$.

\newcommand{\VV}[2]{\mathcal{V}^{#1}_{#2}}
\newcommand{\HH}[2]{\mathcal{H}^{#1}_{#2}}
We use the following compact notations:
$$\VV{j}{m}:=\Phi^{v\Pq{\mu^{j}_{m\dt}}}_{\dt/2}\qquad \HH{j}{m}:=h\Pq{\mu^{j}_{m\dt}}.$$

We have
\bqn
\gw{\mu^{k+1}_{(n+1)\dt},\mu^k_{(n+1)\dt}}&=&\eqnn
\bqn
&=&
\gw{\VV{k+1}{n+\frac12}\#\Pt{\VV{k+1}{n}\#\mu^{k+1}_{n\dt}+\dtd\HH{k+1}{n}}+\dtd\HH{k+1}{n+\frac12},\VV{k}{n}\#\VV{k}{n}\#\mu^k_{n\dt}+\dtd\HH{k}{n}+\dtd\HH{k}{n}}\leq\nn
&\leq&\gw{\VV{k+1}{n+\frac12}\#\VV{k+1}{n}\#\mu^{k+1}_{n\dt},\VV{k}{n}\#\VV{k}{n}\#\mu^k_{n\dt}}+\dtd\gw{\VV{k+1}{n+\frac12}\# \HH{k+1}{n},\HH{k}{n}}+\dtd\gw{\HH{k+1}{n},\HH{k}{n}}.
\eqnn
We estimate the three parts using Proposition \ref{p-utile}. Estimates will be given for sufficiently large $k$, i.e. for sufficiently small $\dt$. Moreover, to simplify the notation, we write estimates that hold for all $p\geq 1$.
We estimate the first term as follows:
\bqn
\gw{\VV{k+1}{n+\frac12}\#\VV{k+1}{n}\#\mu^{k+1}_{n\dt},\VV{k}{n}\#\VV{k}{n}\#\mu^k_{n\dt}}\leq&&\eqnn
\bqn
&\leq&
(1+2L\dt)\gw{\VV{k+1}{n}\#\mu^{k+1}_{n\dt},\VV{k}{n}\#\mu^k_{n\dt}}+m (1+L\dt)\dt\|v\Pq{\mu^{k+1}_{(n+\frac12)\dt}}-v\Pq{\mu^k_{n\dt}}\|_{C^0}\leq\nn
&\leq&
(1+2L\dt)\Pt{(1+2L\dt)\gw{\mu^{k+1}_{n\dt},\mu^k_{n\dt}}+m\dt(1+L\dt)\|v\Pq{\mu^{k+1}_{n\dt}}-v\Pq{\mu^k_{n\dt}}\|_{C^0}}+\nn
&&+ m\dt (1+L\dt)N\Pt{\dtd (Mm+P)+\gw{\mu^{k+1}_{n\dt},\mu^k_{n\dt}}}\leq\nn
&\leq& (1+(5L+4mN)\dt)\gw{\mu^{k+1}_{n\dt},\mu^k_{n\dt}}+ mN(Mm+P) \dt^2.\eqnn
The second term is estimated as follows, using \r{e-rough}:
\bqn
\gw{\VV{k+1}{n+\frac12}\# \HH{k+1}{n},\HH{k}{n}}\leq \gw{\VV{k+1}{n+\frac12}\# \HH{k+1}{n},\HH{k+1}{n}}+\gw{\HH{k+1}{n},\HH{k}{n}}\leq \dtd M P+Q \gw{\mu^{k+1}_{n\dt},\mu^k_{n\dt}}
\eqnn
The third term is simply estimated by $\gw{\HH{k+1}{n},\HH{k}{n}}\leq Q \gw{\mu^{k+1}_{n\dt},\mu^k_{n\dt}}$. Summing up, we have
\bqn
\gw{\mu^{k+1}_{(n+1)\dt},\mu^k_{(n+1)\dt}}&\leq& (1+C_1\dt)\gw{\mu^{k+1}_{n\dt},\mu^k_{n\dt}}+C_2\dt^2
\eqnn
with $C_1:=5L+4mN+Q$, $C_2:=mN(Mm+P)+\frac{MP}{4}$, both independent on $k,n$. Applying it recursively in $n$ and recalling that $\gw{\mu^{k+1}_0,\mu^k_0}=0$, we have
\bqn
\gw{\mu^{k+1}_{n\dt},\mu^k_{n\dt}}&\leq & C_2\dt^2 \frac{(1+C_1\dt)^n-1}{1+C_1\dt-1}\leq 2n C_2\dt^2.
\eqnl{e-basico}

We use such estimate to prove the convergence of $\mu^k_t$ for each $t\in\Pq{0,1}$. We study the three following cases:
\bd
\i[$\bf t=1$:] it corresponds to $n={\dt}^{-1}$. Applying the triangular inequality and \r{e-basico}, we have
\bqn
\gw{\mu^{k}_1,\mu^{k+l}_1}&\leq &2C_2 \Pt{\frac{1}{2^k}+\frac{1}{2^{k+1}}+\ldots+\frac{1}{2^{k+l}}}\leq \frac{4C_2}{2^k}.
\eqnn
Thus, $\mu^k_1$ is a Cauchy sequence in a complete space, hence it is convergent.
\i[$\bf t=\frac{r}{2^l}$ for some $r,l$ integers:] Consider $\mu^k_t$ starting from $k=l$ and apply estimates similar to the previous case.
\i[any $t\in\Pq{0,1}$:] For each $k$, let $n_k$ be the maximum integer such that $t-n_k\dt\geq 0$. We have
\bqn
\gw{\mu^k_t,\mu^{k+l}_t}&\leq&\gw{\mu^k_t,\mu^k_{n_k2^{-k}}}+\gw{\mu^k_{n_k2^{-k}},\mu^{k+l}_{n_k2^{-k}}}+\gw{\mu^{k+l}_{n_{k+l}2^{-k-l}},\mu^{k+l}_t}\leq\nn
&\leq& 2^{-k}(Mm+P)+ \gw{\mu^k_{n_k2^{-k}},\mu^{k+l}_{n_k2^{-k}}}+2^{-k}(Mm+P).
\eqnn
The first and third term are converging with respect to $k$, uniformly in $l$. The same holds for the second term, as proved at the previous step.
\ed

\subsection{Existence and uniqueness of the solution}

\newcommand{\mt}{\bar \mu_t}

In this section we prove that the measure $\mt:=\lim_k \mu^k_t$ given as a limit of the sequence of the Lagrangian scheme is indeed a weak solution of \r{e-cauchy}. We will then prove that it is unique.

First observe that $\mt$ exists, since it is the limit of a Cauchy sequence in $(\Mu,\gwbase)$, that is complete. It is evident that, since $\mu^k_0=\mu_0$ for all $k$, then $\bar\mu_0=\mu_0$. We have now to prove that, for each $f\in C^\infty_c(\Pt{0,1}\times {\R^d})$, it holds
\bqn
&&\int_0^1 dt \int_{\R^d} \Pt{d\mt\Pt{\partial_t f + v\Pq{\mt}\cdot \nabla f}+d h\Pq{\mt} f}=0
\eqnl{e-weak}

The idea is to use again $\mu^k$. First of all, we recall the solution of the continuity equation with time-independent $v$ or $h$ and with the other part ($h$ or $v$ respectively) being identically zero. For $h\equiv 0$, the solution of \r{e-cauchy} with time-independent $v$ is $\mu_t=\Phi^v_t\# \mu_0$. For $v\equiv 0$, the solution of \r{e-cauchy} with time-independent $h$ is $\mu_t=\mu_0+t\,h$. This gives
\bqn
&&\sum_{n=0}^{2^k-1} \int_{n\dt}^{(n+1)\dt} dt \int_{\R^d} d\mu^k_{t}\Pt{\partial_t f + v\Pq{\mu^k_{n\dt}}\cdot \nabla f}+d h\Pq{\mu^k_{n\dt}} f=\nn
&&\sum_{n=0}^{2^k-1} \int_{n\dt}^{(n+1)\dt} dt \int_{\R^d} \Pt{d \Pt{\Phi^{v\Pq{\mu^k_{n\dt}}}_{t-n\dt}\#\mu^k_{n\dt}+(t-n\dt)h\Pq{\mu^k_{n\dt}}}
\Pt{\partial_t f + v\Pq{\mu^k_{n\dt}}\cdot \nabla f}+d h\Pq{\mu^k_{n\dt}} f}=\eqnn
\bqn
&=&\sum_{n=0}^{2^k-1} \int_{n\dt}^{(n+1)\dt} dt \int_{\R^d} d \Pt{\Phi^{v\Pq{\mu^k_{n\dt}}}_{t-n\dt}\#\mu^k_{n\dt}}\Pt{\partial_t f + v\Pq{\mu^k_{n\dt}}\cdot \nabla f}+\nn
&&+ \sum_{n=0}^{2^k-1} \int_{n\dt}^{(n+1)\dt} dt \Pt{(t-n\dt)\int_{\R^d} dh\Pq{\mu^k_{n\dt}} \partial_t f+\int_{\R^d} dh\Pq{\mu^k_{n\dt}} f}+\nn
&&+ \sum_{n=0}^{2^k-1} \int_{n\dt}^{(n+1)\dt} dt (t-n\dt)\int_{\R^d} dh\Pq{\mu^k_{n\dt}}\,v\Pq{\mu^k_{n\dt}}\cdot \nabla f=\nn
&=&0+0+\sum_{n=0}^{2^k-1} \int_{n\dt}^{(n+1)\dt} dt (t-n\dt)\int_{\R^d} dh\Pq{\mu^k_{n\dt}}\,v\Pq{\mu^k_{n\dt}}\cdot \nabla f.
\eqnn
Due to the boundedness of $h,v,\nabla f$, it exists $C_3$ such that $\Pabs{\int_{\R^d} dh\Pq{\mu^k_{n\dt}}\,v\Pq{\mu^k_{n\dt}}\cdot \nabla f(t,.)}\leq C_3$, that gives 
\bqn
&&\Pabs{\lim_k \sum_{n=0}^{2^k-1} \int_{n\dt}^{(n+1)\dt} dt (t-n\dt)\int_{\R^d} dh\Pq{\mu^k_{n\dt}}\,v\Pq{\mu^k_{n\dt}}\cdot \nabla f}\leq\nn
&&\leq \lim_k C_3 \sum_{n=0}^{2^k-1} \int_{n\dt}^{(n+1)\dt} dt (t-n\dt)=C_3\,2^k \frac{\dt^2}{2}=0
\eqnn

Going back to prove \r{e-weak}, we prove equivalently that
\bqn
&\lim_k& \left|\int_0^1 dt \int_{\R^d} \Pt{d\mt\Pt{\partial_t f + v\Pq{\mt}\cdot \nabla f}+d h\Pq{\mt} f}+\right.\nn
&&- \left.\sum_{n=0}^{2^k-1} \int_{n\dt}^{(n+1)\dt} dt \int_{\R^d} d\mu^k_{t}\Pt{\partial_t f + v\Pq{\mu^k_{n\dt}}\cdot \nabla f}+d h\Pq{\mu^k_{n\dt}} f\right|=0.
\eqnn

Since $\lim_k\gw{\mu^k_t,\mt}=0$ implies $\mu^k_t\weak \mt$ (Theorem \ref{t-convergence}), then $$\lim_k \sum_{n=0}^{2^k-1} \int_{n\dt}^{(n+1)\dt} dt \int_{\R^d} d\mu^k_{t}\partial_t f=\int_0^1 dt \int_{\R^d} d\mt \partial_t f.$$ The same applies to $h$, for which we have $\lim_k\gw{h\Pq{\mu^k_t},h\Pq{\mt}}=0$, thus 
$$\lim_k \sum_{n=0}^{2^k-1} \int_{n\dt}^{(n+1)\dt} dt \int_{\R^d} dh\Pq{\mu^k_{t}}\,f=\int_0^1 dt \int_{\R^d} dh\Pq{\mt} \, f.$$

We are now left to prove $$\lim_k \left|\int_0^1 dt \int_{\R^d} d\mt v\Pq{\mt}\cdot \nabla f-\sum_{n=0}^{2^k-1} \int_{n\dt}^{(n+1)\dt} dt \int_{\R^d} d\mu^k_{t} v\Pq{\mu^k_{n\dt}}\cdot \nabla f\right|=0.$$
Using the triangular inequality, we prove it by proving the three following:
\bi
\i $\lim_k \left|\sum_{n=0}^{2^k-1} \int_{n\dt}^{(n+1)\dt} dt \int_{\R^d} \Pt{d\mt-d\mu^k_{t}} v\Pq{\mt}\cdot \nabla f\right|=0$.
Since $\mt$ is continuous and $v$ is Lipschitz with respect to its argument and bounded, then $v\Pq{\mt}\cdot \nabla f \in C^\infty_c(\Pt{0,1}\times {\R^d})$, thus $\mu^k_t\weak \mt$ implies the result.
\i $\lim_k \left|\sum_{n=0}^{2^k-1} \int_{n\dt}^{(n+1)\dt} dt \int_{\R^d} d\mu^k_{t} \Pt{v\Pq{\mt}-v\Pq{\mu^k_{t}}}\cdot \nabla f\right|=0$.
Since it exists a constant $C_4$ such that $\gw{\mu^k_t,\mt}\leq 2^{-k} C_4$, then the previous limit can be proved by recalling that $\mu^k_t$ has finite mass smaller than $m$ and observing that $\nabla f$ is a bounded function. Thus 
\bqn
&&\lim_k \left|\sum_{n=0}^{2^k-1} \int_{n\dt}^{(n+1)\dt} dt \int_{\R^d} d\mu^k_{t} \Pt{v\Pq{\mt}-v\Pq{\mu^k_{t}}}\cdot \nabla f\right|\leq\lim_k \left| \int_0^1 dt \, m\, 2^{-k} N C_4\, \max_{\Pq{0,1}\times {\R^d}} \Pt{\nabla f}\right|=0.
\eqnn
\i $\lim_k \left|\sum_{n=0}^{2^k-1} \int_{n\dt}^{(n+1)\dt} dt \int_{\R^d} d\mu^k_{t} \Pt{v\Pq{\mu^k_t}-v\Pq{\mu^k_{n\dt}}}\cdot \nabla f\right|=0$.
The proof is similar to the previous case, recalling that $\|v\Pq{\mu^k_t}-v\Pq{\mu^k_{n\dt}}\|_{C_0}\leq (t-n\dt)(Mm+P)$.
\ei

We have thus proved the first part of the following proposition
\bp The measure $\mt:=\lim_k \mu^k_t$ is a weak solution of \r{e-cauchy}. Moreover, $\mt\in\Muu$.\ep
\bproof We have to prove that $\mt\in\Muu$. Define the non-autonomous vector field $v_t:=v\Pq{\mt}$ and the time-dipendent measure $h_t:=h\Pq{\mu_t}$. Recall that $\mt$ is continuous with respect to time, thus $v_t$ is a continuous vector field with respect to time, each $h_t$ is absolutely continuous and $h_t$ is continuous with respect to time. The corresponding unique solution of \r{e-cauchy}  is thus in $\Muu$, see Theorem ref{t-misuretrasporto}. Hence, $\mt\in\Muu$.
\eproof

We now prove that we have continuous dependence of the solution of \r{e-cauchy} from the initial data. The estimate also give uniqueness of the solution of \r{e-cauchy}.

\bt Let $\mu_t,\nu_t$ be two solutions of \r{e-cauchy}, with initial data $\mu_0,\nu_0\in\Muu$ respectively. Let $v,h$ satisfy \Hp. Then
\bqn
\gw{\mu_t,\nu_t}\leq e^{t\Pt{\frac{p+1}{p}L+2mN+Q+1}}\gw{\mu_0,\nu_0}.
\eqnl{e-gronsol}
In particular, under assumption \Hp, the solution of \r{e-cauchy} is unique.
\et
\bproof The key observation is that a solution of \r{e-cauchy} at time $t+s$ is of the form $\mu_{t+s}=\Phi^{v\Pq{\mu_t}}_s\#\mu_t+s\, h\Pq{\mu_t}+o(s)$. Also observe that the mass of each solution satisfies $\mu_t(\R^d)\leq \mu_0(\R^d)+tP\leq m^p,$ with $m:=(\mu_0(\R^d)+P)^{1/p}$.

Take now $v,h$ satisfying \Hp, and $\mu_t,\nu_t$ two solutions of \r{e-cauchy}, with initial data $\mu_0,\nu_0$ respectively (eventually coinciding). Define $\eps(t):=\gw{\mu_t,\nu_t}$ and observe that it is a Lipschitz function, since
\bqn
e(t+s)-e(t)&=&\gw{\mu_{t+s},\nu_{t+s}}-\gw{\mu_t,\nu_t}\leq\nn
&\leq&\gw{\mu_{t+s},\mu_t}+\gw{\mu_t,\nu_t}+\gw{\nu_t,\nu_{t+s}}-\gw{\mu_t,\nu_t}\leq 2s(Mm+P)
\eqnn
Now observe the following
\bqn
\eps(t+s)&\leq& \gw{\Phi^{v\Pq{\mu_t}}_s\#\mu_t,\Phi^{v\Pq{\nu_t}}_s\#\nu_t}+sQ\gw{\mu_t,\nu_t}+o(s)\leq\Pt{1+{\scriptstyle \frac{p+1}{p}}Ls +o(s)}\gw{\mu_t,\nu_t}+\nn
&& +\mu_t(\R^d)^{1/p}(1+Ls/p+o(s))(s+o(s)) N \gw{\mu_t,\nu_t}+ sQ\gw{\mu_t,\nu_t}+o(s)\leq\nn
&\leq& \eps(t)+s\Pt{{\scriptstyle \frac{p+1}{p}}L+2mN+Q}\eps(t)+o(s)\leq \eps(t)+s\Pt{{\scriptstyle \frac{p+1}{p}}L+2mN+Q+1}\eps(t).
\eqnn
The last estimate holds for sufficiently small $s$. Using the integral form of the Gronwall inequality, we have \r{e-gronsol}. The uniqueness of the solution of \r{e-cauchy} is a direct consequence.
\eproof

\noindent \b{Acknowledgments}: This work was conducted during a  visit of F. Rossi to Rutgers University, Camden, NJ, USA. He thanks the institution for its hospitality.

\end{document}